\documentclass[11pt,a4paper]{article}
\usepackage[T1]{fontenc}
\usepackage[utf8]{inputenc}
\usepackage{lmodern}
\usepackage{amsmath,amssymb,mathtools,amsthm}
\usepackage[margin=28mm]{geometry}
\usepackage[numbers,sort&compress]{natbib}
\usepackage{microtype}
\usepackage{xurl}
\usepackage[hidelinks]{hyperref}
\usepackage{doi}
\renewcommand{\doi}[1]{\url{https://doi.org/#1}}
\numberwithin{equation}{section}
\newcommand{\Z}{\mathbb Z}
\newcommand{\Q}{\mathbb Q}
\newcommand{\R}{\mathbb R}
\newcommand{\C}{\mathbb C}
\newcommand{\F}{\mathbb F}
\newcommand{\Zp}{\mathbb Z_{\ge0}}

\newcommand{\Nrm}{\operatorname{N}}

\DeclareMathOperator{\Impart}{Im}
\DeclareMathOperator{\Repart}{Re}

\theoremstyle{plain}
\newtheorem{theorem}{Theorem}[section]
\newtheorem{lemma}[theorem]{Lemma}
\newtheorem{proposition}[theorem]{Proposition}
\newtheorem{corollary}[theorem]{Corollary}
\theoremstyle{definition}
\newtheorem{definition}[theorem]{Definition}
\theoremstyle{remark}

\allowdisplaybreaks[1]
\title{Cubic Optimization over Three Nonnegative Integers Is NP-hard\thanks{%
A version of this manuscript was posted on
\href{https://optimization-online.org/2026/11/cubic-optimization-over-three-nonnegative-integers-is-np-hard/}{Optimization Online}
on November~1, 2026.}}
\author{Linchuan Wei\\[0.4em]
\small School of Data Science\\
\small The Chinese University of Hong Kong, Shenzhen\\
\small Shenzhen, China\\
\small\href{mailto:linchuanwei@cuhk.edu.cn}{\texttt{linchuanwei@cuhk.edu.cn}}}
\date{November 4, 2026}
\hypersetup{
  pdftitle={Cubic Optimization over Three Nonnegative Integers Is NP-hard},
  pdfauthor={Linchuan Wei},
  pdfsubject={Fixed-dimensional nonlinear integer optimization},
  pdfkeywords={Integer polynomial optimization, Fixed dimension, Diophantine equations, NP-hardness}}

\begin{document}
\maketitle

\begin{abstract}
We prove that deciding whether an integer cubic polynomial takes a value
at most one over three nonnegative integers is NP-hard. Hardness holds
even when every feasible objective value is a positive integer, so every
constructed instance has a finite, attained minimum. This establishes a
sharp dimension boundary for the tractability of cubic integer
optimization in fixed dimension.
Our reduction from positive 1-in-3 satisfiability encodes Boolean
constraints as modular conditions and uses an algebraic norm to
construct a cubic objective whose minimum is one exactly for
satisfiable formulas. The resulting six-parameter threshold language is
NP-complete and admits short exponent certificates.
\end{abstract}

\par\medskip

\section{Introduction}\label{sec:intro}
The complexity of integer polynomial optimization depends on both the
degree of the objective and the number of variables.
Lenstra~\cite{Lenstra1983} proved that integer linear optimization over
rational polyhedra is polynomial-time solvable in every fixed dimension.
For quadratic objectives, recent preprints by
Wei~\cite{Wei2026} and by Ari and
Hildebrand~\cite[Theorem 1.1]{AriHildebrand2026} independently establish
polynomial-time solvability in every fixed dimension. Their algorithms
allow indefinite rational objectives and unbounded rational polyhedra.

At degree four, fixed dimension no longer ensures tractability.
De Loera, Hemmecke, K\"oppe, and Weismantel~\cite{DeLoeraEtAl2006}
show that minimizing a quartic polynomial over the integer points of a
two-dimensional rectangle is NP-hard. Cubic objectives therefore lie
between the tractable quadratic case and quartic hardness.

For cubic objectives, Del Pia, Hildebrand, Weismantel, and
Zemmer~\cite{DelPiaEtAl2016} give a polynomial-time minimization
algorithm over the integer points of a rational polyhedron in the plane.
Whether cubic integer optimization over
possibly unbounded rational polyhedra in dimension three is NP-hard
remained open. We resolve this question by proving that integer cubic
optimization over three nonnegative integers is NP-hard. Together with
the polynomial-time algorithms for linear and quadratic objectives in
every fixed dimension and for cubic objectives in the plane, and the
NP-hardness of quartic optimization in the plane, we complete the
tractability/intractability landscape of integer polynomial optimization
in fixed dimension.

\subsection{Problem and main result}
Define the cubic threshold language
\begin{equation}\label{eq:threshold}
\mathcal L_3=\left\{(f,\gamma):
\begin{array}{l}
f\in\Z[x,y,z],\quad \deg f\le3,\quad\gamma\in\Z,\\
\exists(x,y,z)\in\Zp^3:\ f(x,y,z)\le\gamma
\end{array}\right\}.
\end{equation}
Rational data can be converted to integer data by multiplying the inequality by a positive common denominator.

\begin{theorem}\label{thm:main}
The language \(\mathcal L_3\) is NP-hard. This remains true even when
\(f\) has a finite, attained minimum over \(\Zp^3\).
\end{theorem}

\subsection{Relation to previous work}
Ari and Hildebrand~\cite[Theorem 1.2]{AriHildebrand2026} develop a
polynomial-time algorithm for integer cubic optimization over bounded
rational polyhedra in every fixed dimension. Their algorithm constructs
a symmetric displacement cover with polynomially many dyadic reflecting
cells. The cubic identities
provide an admissibility certificate on each cell: every global
minimizer in the cell satisfies the admissibility condition, supporting
inequalities hold among admissible integer points, and negative-curvature
search supplies linear curvature cuts separating inadmissible integer
query points. These supporting inequalities and curvature cuts enable
convex integer optimization methods to minimize the objective over the
admissible integer points in each cell. In contrast, our hardness result
holds on the fixed unbounded orthant \(\Zp^3\).

Earlier arithmetic reductions show the NP-hardness of integer quartic
optimization.
Manders and Adleman~\cite{MandersAdleman1978} prove NP-completeness of
deciding, for given positive integers \(a,b,c\), whether
\(a x^2+b y=c\) has a nonnegative integer solution. They also establish
NP-completeness of deciding whether \(x^2\equiv a\pmod b\) has an
integer solution with \(1\le x<c\). De Loera et
al.~\cite{DeLoeraEtAl2006} derive quartic optimization hardness from
the latter problem by minimizing \((x^2-a-by)^2\) over the integer
points of a suitable bounded rectangle. This converts arithmetic
feasibility into optimization with only linear constraints, but squaring
the residual raises the objective degree to four and therefore does not
establish cubic hardness.

We construct a cubic norm equation whose integer solutions correspond
to powers of a chosen algebraic number. This correspondence follows
from Nagell's theorem on fundamental units~\cite{Nagell1930}; see
Louboutin~\cite{Louboutin2015} for a modern treatment. We give a direct
proof for our parameter family in Section~\ref{sec:norm} and use this
correspondence to encode congruences within a cubic objective.

\subsection{Proof idea}
The reduction starts from positive 1-in-3 satisfiability, in which every
clause contains three distinct Boolean variables and requires
exactly one to be true. We encode all truth values by the remainders of
a single integer at distinct primes. Arithmetic tests enforce
Boolean values and the clause requirements, and the Chinese remainder
theorem combines these tests into one congruence.

A change of variables incorporates this congruence into the auxiliary
cubic equation described above. The resulting cubic objective takes
positive integer values over three nonnegative integers and attains a
minimum of one if and only if the congruence has a solution.

Section~\ref{sec:norm} develops the auxiliary cubic equation and its
integer solutions. Section~\ref{sec:filters} constructs the arithmetic
tests. Section~\ref{sec:reduction} combines them and proves
Theorem~\ref{thm:main}, including the running-time bound.
Corollary~\ref{thm:NPcomplete} records NP-completeness for the
six-parameter family using short exponent certificates.

\subsection{Notation and conventions}\label{subsec:notation}
We write \(\Z\) for the integers, \(\Q\) for the rational numbers,
and \(\Zp=\{0,1,2,\ldots\}\) for the nonnegative integers.
For each positive integer \(d\), \(\Zp^d\) is the set of
\(d\)-tuples of nonnegative integers, and \(\Z_{>0}\) denotes
the positive integers.

The symbols \(\Z[T]\) and \(\Q[T]\) denote polynomial rings in
\(T\) with integer and rational coefficients, respectively. We use
the same convention for several variables, as in \(\Z[x,y,z]\).
For an algebraic number \(\alpha\), \(\Z[\alpha]\) is the ring
of values \(f(\alpha)\) with \(f\in\Z[T]\).

\section{Encoding a congruence by a cubic objective}\label{sec:norm}
The aim of this section is to construct a cubic polynomial whose minimum
over \(\Zp^3\) is one exactly when an associated congruence has a solution.

Call an integer tuple \((A,B,a,b,r,M)\) admissible when
\begin{equation}\label{eq:admissible-star}
M\ge2,\qquad 0\le B,a,b,r<M,\qquad A\ge(M+1)^2.
\end{equation}
Fix an admissible tuple throughout this section.

\subsection{Integer coordinates and the norm form}

We start with the polynomial
\begin{equation}\label{eq:p}
p(T)=T^3-AT^2-BT-1.
\end{equation}

\begin{lemma}\label{lem:cubic-ring}
For every admissible tuple, \(p\) is irreducible over \(\Q\) and has
a real root in \((A,A+1)\). For any such root \(\alpha\), we have
\(\alpha>9\), and every element of \(R=\Z[\alpha]\) has a
unique expression \(U+V\alpha+W\alpha^2\), with \(U,V,W\in\Z\).
\end{lemma}
\begin{proof}
\emph{Irreducibility and a real root.}
To prove that \(p\) is irreducible over \(\Q\), we first rule out
rational roots. Suppose that \(p\) has a rational root. After
canceling common factors and choosing a positive denominator, write
it as \(c/d\), with \(c\in\Z\), \(d\in\Z_{>0}\), and
\(\gcd(c,d)=1\). This last condition means that \(c\) and \(d\)
have no common positive divisor other than \(1\). Clearing denominators
gives
\[
c^3=d(Ac^2+Bcd+d^2).
\]
Thus \(d\mid c^3\), and coprimality forces \(d=1\). The root
must therefore be an integer \(c\). Now \(p(c)=0\) gives
\(c(c^2-Ac-B)=1\), so \(c\) divides \(1\) and hence
\(c\in\{1,-1\}\). Admissibility gives \(A>B\ge0\), which
excludes both candidates:
\[
p(1)=-A-B<0,\qquad p(-1)=B-A-2<0.
\]
A reducible cubic over \(\Q\) must have a linear factor and thus
a rational root. Hence \(p\) is irreducible.

Moreover, \(p(A)=-AB-1<0\) and
\(p(A+1)=(A+1)(A+1-B)-1>0\), so continuity gives a real root
in \((A,A+1)\). Let \(\alpha\) be any root in this interval.
Since \(A\ge9\), we have \(\alpha>9\).

\emph{Integer coordinates.} Since \(p\) is irreducible and
\(p(\alpha)=0\), it is the minimal polynomial of \(\alpha\) over
\(\Q\), of degree three. Every element of \(R\) is
\(f(\alpha)\) for some \(f\in\Z[T]\). The relation
\[
\alpha^3=A\alpha^2+B\alpha+1
\]
reduces every higher power of \(\alpha\) to an integer combination
of \(1,\alpha,\alpha^2\). Equivalently, division of \(f\) by the
monic polynomial \(p\) preserves integer coefficients and gives a
remainder \(U+VT+WT^2\), with \(U,V,W\in\Z\). Since
\(p(\alpha)=0\), evaluating this remainder at \(\alpha\) gives
the required expression for \(f(\alpha)\).

For uniqueness, subtracting two such expressions gives a polynomial
of degree at most two that vanishes at \(\alpha\). The minimal
polynomial has degree three, so the difference polynomial is zero.
Thus the three integer coefficients agree term by term.
\end{proof}

Fix such a root \(\alpha\) and set \(R=\Z[\alpha]\) for the rest
of the section. The integer coordinates supplied by
Lemma~\ref{lem:cubic-ring} will be the variables of the norm form.

The images of \(1,\alpha,\alpha^2\) under multiplication by \(\alpha\)
are \(\alpha\), \(\alpha^2\), and
\(1+B\alpha+A\alpha^2\). Their coordinate columns give the matrix
\begin{equation}\label{eq:matrix}
C=\begin{pmatrix}0&0&1\\1&0&B\\0&1&A\end{pmatrix}.
\end{equation}

\begin{definition}[Norm on the cubic ring]\label{def:norm}
Let \(u=U+V\alpha+W\alpha^2\in R\), where \(U,V,W\in\Z\).
Multiplication by \(u\) is a \(\Z\)-linear map \(R\to R\) with matrix
\(UI+VC+WC^2\), where \(I\) is the \(3\times3\) identity matrix.
The norm of \(u\) is the determinant of this map:
\begin{equation}\label{eq:normdet}
\Nrm(U+V\alpha+W\alpha^2)=N_{A,B}(U,V,W):=\det(UI+VC+WC^2).
\end{equation}
\end{definition}

Explicitly, the matrix of multiplication by \(u\) is
\[
\begin{pmatrix}
U&W&V+AW\\
V&U+BW&BV+(1+AB)W\\
W&V+AW&U+AV+(A^2+B)W
\end{pmatrix}.
\]
Expanding its determinant gives
\begin{equation}\label{eq:normpoly}
\begin{aligned}
N_{A,B}(U,V,W)={}&U^3+V^3+W^3+AU^2V+(A^2+2B)U^2W\\
&-BUV^2-(AB+3)UVW+(B^2-2A)UW^2\\
&+AV^2W-BVW^2.
\end{aligned}
\end{equation}
Multiplication by a product is the composition of the two multiplication
maps. Taking determinants shows that
\(\Nrm(uv)=\Nrm(u)\Nrm(v)\). In particular,
\begin{equation}\label{eq:unit}
\Nrm(\alpha)=\det C=1,\qquad
\alpha^{-1}=\alpha^2-A\alpha-B\in R.
\end{equation}

\subsection{Classifying positive elements of norm one}

The following classification is a special case of Nagell's theorem for
cubic orders of negative discriminant~\cite[Satz~XXII]{Nagell1930};
see also~\cite[Theorem~8]{Louboutin2015}. We give a direct proof for
admissible tuples. The classification will let us replace the norm-one
equation by a condition on the coefficients of powers of \(\alpha\).

\begin{lemma}\label{lem:units}
For every admissible tuple,
\begin{equation}\label{eq:allunits}
\{u\in R:u>0,\ \Nrm(u)=1\}=\{\alpha^k:k\in\Z\}.
\end{equation}
\end{lemma}
\begin{proof}
By~\eqref{eq:unit}, every integer power of \(\alpha\), including a
negative power, belongs to the left side of~\eqref{eq:allunits}.
For the reverse inclusion, we normalize an element of norm one into
\([1,\alpha)\) and use bounds at the complex embeddings to show that
the normalized element must be one.

\emph{Conjugate bounds.} Set \(\delta=\alpha-A\). The equation \(p(\alpha)=0\) gives
\[
\delta=\frac B\alpha+\frac1{\alpha^2}>0,\qquad
\delta\sqrt\alpha=\frac{B+1/\alpha}{\sqrt\alpha}
<\frac M{\sqrt A}\le\frac M{M+1}<1.
\]
Here \(B\le M-1\) and \(\alpha>A\ge(M+1)^2>1\) give
\(B+1/\alpha<M\). Since \(\delta>0\), squaring the bound
\(\delta\sqrt\alpha<1\) gives \(\delta^2<1/\alpha\).
The factorization
\[
p(T)=(T-\alpha)(T^2+\delta T+1/\alpha)
\]
shows that the other roots of \(p\) solve
\(T^2+\delta T+1/\alpha=0\). This quadratic has negative
discriminant, since
\[
\delta^2-\frac4\alpha
<\frac1\alpha-\frac4\alpha=-\frac3\alpha<0.
\]
Its roots therefore form a nonreal conjugate pair
\(\beta,\overline\beta\), and the quadratic formula gives them as
\[
-\frac\delta2\pm i\sqrt{\frac1\alpha-\frac{\delta^2}{4}}.
\]
Consequently,
\[
|\Impart\beta|^2
=\frac1\alpha-\frac{\delta^2}{4}
>\frac1\alpha-\frac1{4\alpha}
=\frac3{4\alpha}.
\]
Taking positive square roots gives
\(|\Impart\beta|>\sqrt3/(2\sqrt\alpha)>1/(2\sqrt\alpha)\),
the bound we will use below. The same formula gives
\(\Repart\beta=-\delta/2\) and \(|\beta|^2=1/\alpha\).
The characteristic polynomial of \(C\) is \(p\), with three distinct roots.
Thus the eigenvalues of \(UI+VC+WC^2\) are the values of
\(U+VT+WT^2\) at \(\alpha,\beta,\overline\beta\).
Taking their product gives, for
\(u=U+V\alpha+W\alpha^2\) and \(u'=U+V\beta+W\beta^2\),
\begin{equation}\label{eq:normproduct}
\Nrm(u)=u\,u'\overline{u'}=u|u'|^2.
\end{equation}

\emph{Normalization and coefficient bounds.} Let \(u\in R\) satisfy
\(u>0\) and \(\Nrm(u)=1\). Since \(\alpha>1\), we can choose
\(k\in\Z\) so that \(v=u\alpha^{-k}\in[1,\alpha)\).
Equation~\eqref{eq:unit} and multiplicativity give \(v\in R\) and
\(\Nrm(v)=1\). By~\eqref{eq:normproduct}, its value \(v'\) under
the embedding \(\alpha\mapsto\beta\) satisfies
\(|v'|=v^{-1/2}\le1\).
Write \(v=c_0+c_1\alpha+c_2\alpha^2\), with \(c_i\in\Z\).
The relation \(\beta^2=-\delta\beta-1/\alpha\) lets us write this
complex value using two real coefficients:
\[
d_0=c_0-c_2/\alpha,\qquad d_1=c_1-\delta c_2,
\qquad v'=d_0+d_1\beta.
\]
Since \(d_0,d_1\) are real, the imaginary part of \(v'\) is
\(d_1\Impart\beta\). The bounds on \(|v'|\) and
\(|\Impart\beta|\) therefore give \(|d_1|<2\sqrt\alpha\).
To bound \(d_0\), take real parts of \(v'=d_0+d_1\beta\)
and use \(\Repart\beta=-\delta/2\):
\[
d_0=\Repart v'-d_1\Repart\beta
=\Repart v'+\frac{\delta d_1}{2}.
\]
Since \(|\Repart v'|\le|v'|\le1\), the triangle inequality gives
\[
\begin{aligned}
|d_0|&\le|\Repart v'|+\frac{\delta|d_1|}{2}\\
&\le1+\frac{\delta|d_1|}{2}
<1+\delta\sqrt\alpha<2.
\end{aligned}
\]
Here we used \(\delta>0\), \(|d_1|<2\sqrt\alpha\), and
\(\delta\sqrt\alpha<1\) in the last two inequalities.
Substituting the definitions of \(d_0,d_1\) into the real value of
\(v\) gives
\[
v=d_0+d_1\alpha+c_2(\alpha^2+\delta\alpha+1/\alpha).
\]
The factor multiplying \(c_2\) exceeds \(\alpha^2\). The triangle
inequality, \(v<\alpha\), and the bounds on \(d_0,d_1\) therefore give
\[
\begin{aligned}
|c_2|
&=\frac{|v-d_0-d_1\alpha|}{\alpha^2+\delta\alpha+1/\alpha}\\
&<\frac{\alpha+2+2\alpha\sqrt\alpha}{\alpha^2}
<\frac3{\sqrt\alpha}<1.
\end{aligned}
\]
For the bound by \(3/\sqrt\alpha\), \(\alpha>9\) gives
\(\alpha+2<3\alpha<\alpha\sqrt\alpha\); the last uses the same bound.
Because \(c_2\) is an integer, it must be zero. We then have
\(c_0=d_0\in\{-1,0,1\}\) and
\(1\le c_0+c_1\alpha<\alpha\). If \(c_1\le-1\), this value is
negative; if \(c_1\ge2\), it exceeds \(\alpha\). For \(c_1=0\),
the only possibility is \(c_0=1\); for \(c_1=1\), it is
\(c_0=-1\). Thus \(v\) is either \(1\) or \(\alpha-1\).
Finally, \(\Nrm(\alpha-1)=\det(C-I)=-p(1)=A+B>1\), so the
second possibility is excluded. Hence \(v=1\) and \(u=\alpha^k\).
\end{proof}

\subsection{A cubic integer program encoding a congruence}

For \(k\ge0\), define the integer coefficients
\begin{equation}\label{eq:coefficients}
\alpha^k=U_k+V_k\alpha+W_k\alpha^2.
\end{equation}
Multiplication by \(\alpha\) gives
\begin{equation}\label{eq:recurrence}
U_{k+1}=W_k,\qquad V_{k+1}=U_k+BW_k,\qquad
W_{k+1}=V_k+AW_k,
\end{equation}
starting at \((U_0,V_0,W_0)=(1,0,0)\). Since \(A,B\ge0\), the
recurrence preserves nonnegativity, so \(U_k,V_k,W_k\ge0\) for
every \(k\ge0\).

\begin{samepage}
We call a linear congruence among \(U_k,V_k,W_k\), such as
\(W_k-aU_k-bV_k\equiv r\pmod M\), a \emph{coefficient congruence}.
Here these are the coefficients of \(\alpha^k\) in the basis
\(1,\alpha,\alpha^2\).
\end{samepage}

Define the objective
\begin{equation}\label{eq:G-star}
G_*(x,y,z)=N_{A,B}\bigl(x,y,ax+by+r+M(z+1)\bigr).
\end{equation}

The next proposition connects the cubic objective to the coefficient
congruence used in the reduction.

\begin{proposition}\label{prop:equivalence}
For every admissible tuple, \(G_*\) is an integer polynomial of degree
exactly three, is at least one throughout \(\Zp^3\), and has a finite,
attained minimum.
Moreover,
\begin{equation}\label{eq:equivalence}
\begin{aligned}
&\exists(x,y,z)\in\Zp^3:\ G_*(x,y,z)\le1\\
&\qquad\Longleftrightarrow\quad
\exists k\ge0:\ W_k-aU_k-bV_k\equiv r\pmod M.
\end{aligned}
\end{equation}
Call an integer exponent \(k\ge0\) accepted if it satisfies the
congruence on the right of~\eqref{eq:equivalence}. Every accepted
exponent can be replaced by an accepted exponent \(3\le k'<M^3+3\)
yielding a nonnegative integer point of value one.
\end{proposition}
\begin{proof}
\emph{Degree and minimum.}
Equations~\eqref{eq:normpoly} and~\eqref{eq:G-star} show that
\(G_*\) has integer coefficients and degree at most three.
The coefficient of \(z^3\) in~\eqref{eq:G-star} is \(M^3\), so the degree is
exactly three. At an integer point in the orthant, the substituted
coefficients satisfy \(U,V\ge0\) and \(W\ge M>0\). Thus
\(u=U+V\alpha+W\alpha^2>1\). Its complex embedding is nonzero: otherwise
the nonzero rational polynomial \(U+VT+WT^2\) would vanish at \(\beta\),
whose minimal polynomial is the irreducible cubic \(p\). By
\eqref{eq:normproduct}, \(G_*>0\), and integrality gives \(G_*\ge1\).
The set of values is a nonempty subset of the positive integers, so it
has a least element, attained at an integer point.

\emph{From a point to a congruence.} If \(G_*\le1\), positivity
forces \(G_*=1\), so the associated element \(u\) has norm one.
Lemma~\ref{lem:units} gives
\(u=\alpha^k\) for \(k\in\Z\); since \(u>1\), actually \(k>0\).
Uniqueness of the coefficient expression identifies
\(U_k=x\), \(V_k=y\), and
\(W_k=ax+by+r+M(z+1)\), which gives the required congruence.

\emph{From a congruence to a point.} For an accepted exponent
\(k\ge0\), invert the substitution in~\eqref{eq:G-star} to obtain
the candidate coordinates
\begin{equation}\label{eq:inverse}
x=U_k,\qquad y=V_k,\qquad
z=\frac{W_k-aU_k-bV_k-r}{M}-1.
\end{equation}
The recurrence gives \(x,y\ge0\), and the congruence makes \(z\)
an integer. Since \(M>0\), the remaining requirement \(z\ge0\)
is equivalent to \(W_k-aU_k-bV_k\ge r+M\). We first prove this
bound for \(k\ge3\), then use periodicity modulo \(M\) to reach
that range from any accepted exponent.

Write \(L_k=W_k-aU_k-bV_k\). The recurrence gives
\(0\le U_j\le W_j\) for every \(j\ge1\), since
\(W_j-U_j=V_{j-1}+(A-1)W_{j-1}\ge0\).
Also \(W_2=1\) and \(W_{j+1}\ge AW_j\), so \(W_j\ge1\) for \(j\ge2\).
Hence, for \(k\ge3\),
\begin{equation}\label{eq:liftbounds}
\begin{aligned}
L_k
&=V_{k-1}-bU_{k-1}+(A-a-bB)W_{k-1}\\
&\ge\bigl(A-a-b(B+1)\bigr)W_{k-1}\\
&\ge A-(M^2-1)\ge2M+2>r+M.
\end{aligned}
\end{equation}
For the first inequality we used \(V_{k-1}\ge0\) and
\(U_{k-1}\le W_{k-1}\). For the second, admissibility gives
\(a+b(B+1)\le(M-1)+(M-1)M=M^2-1\), so the coefficient of
\(W_{k-1}\) is positive, and \(W_{k-1}\ge1\).

The residue triples \((U_k,V_k,W_k)\bmod M\) lie in a set of \(M^3\)
possible triples, since each coordinate has \(M\)
possible residues. By~\eqref{eq:recurrence}, the next coefficient
column is obtained by multiplying the current column by \(C\).
Since \(\det C=1\), the inverse formula divides integer cofactors
by one, so \(C^{-1}\) has integer entries. Explicitly,
\[
C^{-1}=\begin{pmatrix}-B&1&0\\-A&0&1\\1&0&0\end{pmatrix}.
\]
The identities \(CC^{-1}=C^{-1}C=I\) remain valid after reducing
all entries modulo \(M\). Thus multiplication by \(C^{-1}\)
recovers a unique predecessor for every residue triple.
Multiplication by \(C\) is therefore a bijection, or equivalently
a permutation, of these \(M^3\) triples.

Let \(e=(1,0,0)^{\mathsf T}\) be the initial coefficient column.
The column at index \(k\) is \(C^ke\). Among the \(M^3+1\)
columns at indices \(0,\ldots,M^3\), two coincide modulo \(M\),
so for some \(0\le i<j\le M^3\) we have
\[
C^ie\equiv C^je\pmod M.
\]
Multiplying both sides by \(C^{-i}\) undoes the first \(i\)
updates and gives
\[
e=C^{-i}C^ie\equiv C^{-i}C^je=C^{j-i}e\pmod M.
\]
Thus the column at index \(j-i\) equals the initial column modulo
\(M\), and \(1\le j-i\le M^3\). Choose \(h\) to be the least
positive return time, so \(1\le h\le M^3\). For every \(k\ge0\),
multiplying \(C^he\equiv e\pmod M\) by \(C^k\) gives
\[
C^{k+h}e\equiv C^ke\pmod M.
\]
Hence the coefficient sequence modulo \(M\) repeats from index zero
after \(h\) steps, so acceptance depends only on \(k\bmod h\).
For any accepted \(k\ge0\), choose the representative
\(k'\equiv k\pmod h\) in \(\{3,\ldots,h+2\}\).
This interval contains \(h\) consecutive integers, so it represents
every residue class, including those of accepted exponents \(0,1,2\).
Periodicity preserves the congruence, and \(k'\ge3\) puts the new
exponent in the range of~\eqref{eq:liftbounds}.
Apply~\eqref{eq:inverse} with \(k'\) in place of \(k\).
The bound~\eqref{eq:liftbounds} now ensures \(z\ge0\), so this
gives a nonnegative integer point. Substitution recovers the coefficient triple of
\(\alpha^{k'}\), so
\(G_*(x,y,z)=\Nrm(\alpha^{k'})=1\).
Since \(h\le M^3\), the chosen exponent satisfies
\(3\le k'\le h+2<M^3+3\).
\end{proof}

\section{Encoding Boolean constraints by coefficient congruences}\label{sec:filters}
Proposition~\ref{prop:equivalence} reduces the existence of a nonnegative
integer point with \(G_*=1\) to the existence of an exponent \(k\ge0\)
satisfying a coefficient congruence. We now construct tests that encode
a positive 1-in-3 Boolean formula: a conjunction of clauses, each
requiring exactly one of three distinct, unnegated variables to be true.
Variable primes \(q_i\)
encode truth values through \(k\bmod q_i\in\{0,1\}\). A separate field
prime \(\pi\) is assigned to each variable or clause test, which is
constructed in \(\F_\pi=\Z/\pi\Z\). We first convert finite-field power
tests into coefficient congruences, then design tests enforcing Boolean
values and exactly-one clauses. Section~\ref{sec:reduction} will choose
simultaneous integer parameters \(A,B,a,b,r\) and combine all tests into
one coefficient congruence.

\subsection{Converting power tests into coefficient congruences}
Let \(\pi\) be prime. The order of a nonzero element of \(\F_\pi\) is the least positive integer \(h\) for which its \(h\)-th power is one. Let \(g\in\F_\pi^\times\) have order greater than five. Set
\begin{equation}\label{eq:roots}
\lambda_1=g,\qquad \lambda_2=g^2,\qquad \lambda_3=g^{-3}.
\end{equation}
These elements are distinct: an equality would give \(g=1\), \(g^4=1\), or \(g^5=1\). Their product is one. Define residues
\begin{equation}\label{eq:ABresidues}
A_\pi=\lambda_1+\lambda_2+\lambda_3,\qquad
B_\pi=-(\lambda_1\lambda_2+\lambda_1\lambda_3+\lambda_2\lambda_3).
\end{equation}
Whenever \(A\equiv A_\pi\) and \(B\equiv B_\pi\pmod\pi\),
\[
p(T)\equiv(T-\lambda_1)(T-\lambda_2)(T-\lambda_3)\pmod\pi.
\]
Evaluating the polynomial remainder of \(T^k\) at these roots gives
\begin{equation}\label{eq:evaluate}
U_k+V_k\lambda_i+W_k\lambda_i^2=\lambda_i^k\quad\text{in }\F_\pi,
\qquad i=1,2,3.
\end{equation}

Choose weights \(\eta_1,\eta_2,\eta_3\) and a target \(\rho\) in \(\F_\pi\), and put
\begin{equation}\label{eq:ell}
\ell_j=\sum_{i=1}^3\eta_i\lambda_i^j\qquad(j=0,1,2).
\end{equation}
\begin{definition}[Finite-field power test]\label{def:powertest}
For the chosen \(g\), weights \(\eta_1,\eta_2,\eta_3\), and target \(\rho\),
the \emph{finite-field power test} accepts an integer exponent \(k\ge0\) if
\begin{equation}\label{eq:exponentialtest}
\eta_1g^k+\eta_2g^{2k}+\eta_3g^{-3k}=\rho\quad\text{in }\F_\pi.
\end{equation}
\end{definition}
By~\eqref{eq:evaluate}, the test accepts \(k\) exactly when
\(\ell_0U_k+\ell_1V_k+\ell_2W_k=\rho\). If \(\ell_2\ne0\), this condition can be written as
\begin{equation}\label{eq:localcongruence}
W_k-a_\pi U_k-b_\pi V_k=r_\pi\quad\text{in }\F_\pi,
\end{equation}
where
\begin{equation}\label{eq:normalize}
a_\pi=-\ell_0/\ell_2,\qquad b_\pi=-\ell_1/\ell_2,\qquad r_\pi=\rho/\ell_2.
\end{equation}
Thus we only need to design power tests of the form~\eqref{eq:exponentialtest} for which the coefficient \(\ell_2\) of \(W_k\) is nonzero.

\subsection{A test for one Boolean variable}
\begin{lemma}\label{lem:variable}
Let \(q\ge37\) be prime and let \(\pi\equiv1\pmod q\) be prime.
There exist residues \(A_\pi,B_\pi,a_\pi,b_\pi,r_\pi\in\F_\pi\)
with the following property. For any integers \(A,B\) satisfying
\(A\equiv A_\pi\pmod\pi\) and \(B\equiv B_\pi\pmod\pi\), let
\(U_k+V_kT+W_kT^2\) be the remainder of \(T^k\) upon division by
the polynomial \(p(T)\) in~\eqref{eq:p}. Then, for every \(k\ge0\),
\[
W_k-a_\pi U_k-b_\pi V_k=r_\pi\text{ in }\F_\pi
\quad\Longleftrightarrow\quad k\bmod q\in\{0,1\}.
\]
The residues can be found in time polynomial in the numerical value of \(\pi\).
\end{lemma}
\begin{proof}
We construct a power test that accepts exactly the residues zero and one
modulo \(q\), then convert it into the required coefficient congruence.

\emph{Choosing \(g\).}
The multiplicative group \(\F_\pi^\times\) is cyclic of order \(\pi-1\).
Since \(q\mid\pi-1\), a generator \(\zeta\) gives an element
\(h=\zeta^{(\pi-1)/q}\). For any integer \(s\), the condition \(h^s=1\)
is equivalent to \(q\mid s\), because \(\zeta\) has order \(\pi-1\).
Thus \(h\) generates a subgroup of order \(q\). Every element of this
subgroup has its \(q\)-th power equal to one, so its order divides \(q\).
Since \(q\) is prime, each of the \(q-1\) nonidentity elements has order
exactly \(q\). Choose one, denoted by \(g\), for which
\(g^2-g-1\ne0\). Such a choice exists because the quadratic polynomial
\(T^2-T-1\) has at most two roots in the field, leaving at least
\(q-3>0\) eligible elements.

\emph{The power test.}
Since \(g\) has order \(q>5\), the roots in~\eqref{eq:roots} are distinct.
Define \(A_\pi,B_\pi\) by~\eqref{eq:ABresidues}. For any integers \(A,B\)
with these residues, the evaluation identities~\eqref{eq:evaluate} apply
to the remainder coefficients \(U_k,V_k,W_k\).

Set \(t=g^k\). To accept exactly \(t=1\) and \(t=g\), use the equation
\[
(t-1)(t-g)=t^2-(1+g)t+g=0.
\]
This is the power test~\eqref{eq:exponentialtest} with
\[
(\eta_1,\eta_2,\eta_3)=(-(1+g),1,0),\qquad \rho=-g.
\]
A field has no zero divisors, so the factorization shows that these are
the only two accepted values of \(t\). Since \(g\) has order \(q\),
\(g^k=1\) is equivalent to \(q\mid k\), and \(g^k=g\) is equivalent to
\(q\mid(k-1)\). The test therefore accepts precisely
\(k\equiv0,1\pmod q\).

\emph{The coefficient congruence.}
The weights in~\eqref{eq:ell} are
\[
\begin{aligned}
\ell_0&=-(1+g)+1=-g,\\
\ell_1&=-(1+g)g+g^2=-g,\\
\ell_2&=-(1+g)g^2+g^4=g^2(g^2-g-1).
\end{aligned}
\]
Both factors in \(\ell_2=g^2(g^2-g-1)\) are nonzero by our choice of \(g\),
so \(\ell_2\) has an inverse in \(\F_\pi\). By~\eqref{eq:evaluate}, the
power test is equivalent to
\[
-gU_k-gV_k+\ell_2W_k=-g.
\]
Define the normalized residues as in~\eqref{eq:normalize}:
\[
a_\pi=b_\pi=g/\ell_2,\qquad r_\pi=-g/\ell_2.
\]
Dividing the preceding coefficient equation by \(\ell_2\) now gives
\(W_k-a_\pi U_k-b_\pi V_k=r_\pi\). This congruence holds if and only if
the power test accepts \(k\), hence if and only if
\(k\bmod q\in\{0,1\}\).

\emph{Computing the residues.}
To compute these residues, enumerate the nonzero field elements
\(c=1,\ldots,\pi-1\) and test
\[
c\ne1,\qquad c^q=1,\qquad c^2-c-1\ne0.
\]
The first two conditions give order exactly \(q\): the order divides the
prime \(q\) and is not one. Our existence argument guarantees that at
least one candidate also satisfies the third condition. Take the first
such candidate as \(g\).

There are at most \(\pi-1\) candidates. To compute \(c^q\), write
\(q=\sum_{j=0}^{m}b_j2^j\), where \(b_j\in\{0,1\}\) and
\(m=\lfloor\log_2q\rfloor\). Repeated squaring modulo \(\pi\) produces
\(c,c^2,c^4,\ldots,c^{2^m}\) in \(m\) modular multiplications.
Multiplying the powers for which \(b_j=1\) takes at most another \(m\)
multiplications and gives \(c^q\). Thus each power computation uses
\(O(\log q)\) modular multiplications. Reduction modulo \(\pi\) after
each operation keeps the stored residues below \(\pi\). Since \(q<\pi\),
each candidate test takes time polynomial in \(\log\pi\).

Once \(g\) is found, the formulas for \(A_\pi,B_\pi,a_\pi,b_\pi,r_\pi\)
require only a fixed number of field operations, with inverses computed
by the extended Euclidean algorithm. Combining this cost with the search
over at most \(\pi-1\) candidates proves that the construction is
deterministic and takes time polynomial in the numerical value of \(\pi\).
\end{proof}

\subsection{A test for an exactly-one clause}
An exactly-one clause accepts a Boolean assignment when precisely one of
its three variables is true. Fix distinct primes \(q_1,q_2,q_3\ge37\) and a prime
\(\pi\equiv1\pmod{q_1q_2q_3}\). Let
\[
S_i=\{t\in\F_\pi^\times:t^{q_i}=1,\ t\ne1\}\qquad(i=1,2,3).
\]
By cyclicity, \(S_i\) has \(q_i-1\) elements, each of exact order \(q_i\).
For \((u,v,w)\in S_1\times S_2\times S_3\), the product \(g=uvw\) has
order \(q_1q_2q_3\). Indeed, raising \(g\) to this product gives one.
Conversely, if \(g^h=1\), raising to \(q_2q_3\) gives
\(u^{hq_2q_3}=1\), since the powers of \(v,w\) become one.
Thus \(q_1\mid hq_2q_3\), and coprimality gives \(q_1\mid h\).
The other two divisibilities follow in the same way.

To build the clause test, we encode the three Boolean values in the
single field element \(g^k\). Suppose that
\(k\bmod q_i=\xi_i\in\{0,1\}\) for \(i=1,2,3\). Since \(u,v,w\)
have orders \(q_1,q_2,q_3\), respectively,
\[
g^k=u^kv^kw^k=u^{\xi_1}v^{\xi_2}w^{\xi_3}.
\]
Each bit determines whether its corresponding factor appears. Thus the
eight Boolean assignments give the products
\begin{equation}\label{eq:eight}
1,\ u,\ v,\ w,\ uv,\ uw,\ vw,\ uvw.
\end{equation}
For example, the assignment \((1,1,0)\) gives \(uv\).
These products are distinct. If two of them were equal, dividing one by the other
and raising the equality to \(q_2q_3\) would show, as above, that the
difference of their exponents of \(u\) is divisible by \(q_1\). This difference
lies in \(\{-1,0,1\}\), so it is zero; the same argument applies to
\(v,w\). Consequently, an assignment satisfies the clause precisely when
\(g^k\) is \(u,v\), or \(w\). The clause test must therefore
accept these three values and reject the other five in~\eqref{eq:eight}.

\begin{lemma}\label{lem:clause}
There exist residues
\(A_\pi,B_\pi,a_\pi,b_\pi,r_\pi\in\F_\pi\) with the following property.
For any integers \(A,B\) satisfying \(A\equiv A_\pi\pmod\pi\) and
\(B\equiv B_\pi\pmod\pi\), let \(U_k+V_kT+W_kT^2\) be the remainder
of \(T^k\) upon division by the polynomial \(p(T)\) in~\eqref{eq:p}.
Then, for every \(k\ge0\) whose residues \(\xi_i:=k\bmod q_i\) lie in
\(\{0,1\}\) for \(i=1,2,3\),
\[
W_k-a_\pi U_k-b_\pi V_k=r_\pi\text{ in }\F_\pi
\quad\Longleftrightarrow\quad \xi_1+\xi_2+\xi_3=1.
\]
The residues can be found in time polynomial in the numerical
quantity \(\pi+q_1q_2q_3\).
\end{lemma}
\begin{proof}
With \(t=g^k\), multiplying a power test~\eqref{eq:exponentialtest}
by \(t^3\) gives a polynomial containing only the powers
\(t^5,t^4,t^3\) and a constant. We build such a polynomial with factors
\((t-u)(t-v)(t-w)\), then choose
\((u,v,w)\in S_1\times S_2\times S_3\) so that its remaining factor
is nonzero at the other Boolean values and the coefficient test can be normalized.

\emph{Designing the power test.}
For a prospective triple, set
\begin{equation}\label{eq:clauseparameters}
\begin{gathered}
g=uvw,\qquad s_1=u+v+w,\qquad s_2=uv+uw+vw,\\
D=s_2^2-s_1g,\qquad L=s_2g-s_1D,\qquad
\rho=2s_1s_2g-s_2^3-g^2,
\end{gathered}
\end{equation}
and put \(Q(t)=Dt^2+s_2gt+g^2\). These choices give
\begin{equation}\label{eq:clauseidentity}
(t-u)(t-v)(t-w)Q(t)=Dt^5+Lt^4-\rho t^3-g^3.
\end{equation}
Take weights \((\eta_1,\eta_2,\eta_3)=(L,D,-g^3)\) and target \(\rho\).
The resulting power test is
\[
Lg^k+Dg^{2k}-g^3g^{-3k}=\rho.
\]
Since \(t=g^k\ne0\), the test is equivalent to setting the right side
of~\eqref{eq:clauseidentity} equal to zero. Thus the test holds for
\(t=u,v,w\).
At each of the other five Boolean values, the three linear factors are
nonzero by distinctness in~\eqref{eq:eight}. To reject these values,
we therefore require
\begin{equation}\label{eq:avoidQ}
Q(t)\ne0\qquad(t=1,uv,uw,vw,g).
\end{equation}

\emph{Choosing the triple.}
We also need \(\ell_2\ne0\) to convert the power test into a coefficient
congruence. By~\eqref{eq:ell},
\begin{equation}\label{eq:J}
\ell_2=Lg^2+Dg^4-g^{-3}=g^{-3}J,\qquad J=Dg^7+Lg^5-1.
\end{equation}
Because \(g\ne0\), this condition is equivalent to \(J\ne0\).
To satisfy it together with the five conditions~\eqref{eq:avoidQ},
regard \(u,v,w\) as indeterminates and define
\begin{equation}\label{eq:Psi}
\Psi(u,v,w)=JQ(1)Q(uv)Q(uw)Q(vw)Q(uvw).
\end{equation}
At a triple in \(S_1\times S_2\times S_3\), the condition \(\Psi\ne0\)
is equivalent to all six factors being nonzero. We will prove that
\(\Psi\) is a nonzero polynomial whose degree in each variable is
smaller than the size of the corresponding set \(S_i\).

The polynomials \(g,s_1,s_2\) have degree at most one in each variable,
\(D\) has degree at most two, and \(L=s_2g-s_1D\) has degree at most
three. The degree of \(J=Dg^7+Lg^5-1\) in each variable is therefore
at most \(\max\{2+7,3+5\}=9\). Each of the five arguments of \(Q\)
has individual degree at most one, so each substituted \(Q\) has
individual degree at most four. Hence \(\Psi\) has degree at most
\(9+5\cdot4=29\) in each variable.

To prove that \(\Psi\) is a nonzero polynomial, evaluate it at
\(u=v=w=1\). Here
\[
g=1,\quad s_1=s_2=3,\quad D=6,\quad L=-15,\quad
Q(1)=10,\quad J=-10.
\]
Thus \(\Psi(1,1,1)=-10^6\ne0\), because
\(\pi>q_1q_2q_3>5\).

Now \(|S_i|=q_i-1\ge36>29\). Suppose that \(\Psi\) vanished on the
whole grid \(S_1\times S_2\times S_3\). For fixed \(v\in S_2\) and
\(w\in S_3\), the polynomial \(\Psi(u,v,w)\) in \(u\) would have more
roots than its degree, so all its coefficients would be zero. Thus each
coefficient, viewed as a polynomial in \(v,w\), would vanish on
\(S_2\times S_3\). Applying the same root bound first in \(v\) and
then in \(w\) would make every coefficient of \(\Psi\) zero, contrary
to the preceding evaluation. Hence some triple has \(\Psi\ne0\).

\emph{The coefficient congruence.}
For the selected triple, the five conditions~\eqref{eq:avoidQ} hold, so
the power test holds exactly for \(t=u,v,w\) among the eight Boolean values.
Define the roots by~\eqref{eq:roots} and \(A_\pi,B_\pi\) by
\eqref{eq:ABresidues}. Since \(g\) has order \(q_1q_2q_3>5\), the
roots are distinct and the evaluation identities~\eqref{eq:evaluate}
apply. Moreover, \(J\ne0\) gives \(\ell_2\ne0\) by~\eqref{eq:J}.
Define \(a_\pi,b_\pi,r_\pi\) by~\eqref{eq:normalize}. The power test
is then equivalent to \(W_k-a_\pi U_k-b_\pi V_k=r_\pi\), and, for
Boolean residues \(\xi_i\), it accepts precisely when
\(\xi_1+\xi_2+\xi_3=1\).

\emph{Computing the residues.}
Enumerate the nonzero field elements to construct the sets \(S_i\) by
testing \(t^{q_i}=1\) and \(t\ne1\). As in Lemma~\ref{lem:variable},
binary modular exponentiation uses \(O(\log q_i)\) multiplications
per test. Since \(q_i<\pi\), constructing the three sets takes
\(\pi(\log\pi)^{O(1)}\) bit operations. Search at most
\((q_1-1)(q_2-1)(q_3-1)<q_1q_2q_3\) triples for \(\Psi\ne0\). Each
triple requires a fixed number of field operations to evaluate \(\Psi\).
Once a triple is found, compute the five residues from
\eqref{eq:ABresidues} and~\eqref{eq:normalize}, using the extended
Euclidean algorithm for inverses. The denominators \(g\) and \(\ell_2\)
are nonzero, and each field operation takes time polynomial in
\(\log\pi\). The total bit complexity is therefore polynomial in
\(\pi+q_1q_2q_3\).
\end{proof}

\section{The reduction and its running time}\label{sec:reduction}
We reduce from positive 1-in-3 satisfiability, a conjunction of relations
\[
\xi_i+\xi_j+\xi_\ell=1,\qquad \xi_i\in\{0,1\}.
\]
This problem is NP-complete by Schaefer's theorem~\cite[Theorem~2.1]{Schaefer1978}.
We may require three distinct variables in each clause.
To remove repeated occurrences, replace each occurrence by a fresh copy
\(\zeta\) and enforce its equality with the original variable \(\xi\)
using fresh \(s,t\) in
\[
\xi+s+t=1,\qquad \zeta+s+t=1.
\]
Subtraction forces \(\xi=\zeta\). Conversely, when \(\xi=\zeta\),
setting \(s=1-\xi\) and \(t=0\) satisfies both clauses.
Every new clause has distinct variables, and the size increase is
polynomial.

Remove unused variables and relabel the rest. If none remain, add one
unused variable, preserving the satisfiability of the empty conjunction.
Let \(s\ge1\) and \(m\) be the resulting numbers of variables and clauses,
and let \(n\ge2\) be the binary formula length, padded if necessary.
We will use one variable test from Lemma~\ref{lem:variable} for each
variable and one clause test from Lemma~\ref{lem:clause} for each clause,
giving \(K=s+m=O(n)\) local tests.

\subsection{One construction of distinct primes}
We use the following effective form of Linnik's theorem: there are
absolute effective constants \(C_0,L_0>0\) such that for every \(d\ge2\)
and every residue \(c\) coprime to \(d\), some prime \(\pi\) satisfies
\begin{equation}\label{eq:linnik}
\pi\equiv c\pmod d,\qquad \pi\le C_0d^{L_0}.
\end{equation}
This is the general-modulus theorem; see Heath-Brown~\cite[Theorem~6 and the following paragraph]{HeathBrown1992}.
No particular exponent \(L_0\) is needed.

\begin{lemma}\label{lem:primes}
Given odd integers \(d_1,\ldots,d_K\ge3\), where \(K\ge1\), put
\(D=\max_t d_t\). Distinct primes \(\pi_1,\ldots,\pi_K\) satisfying
\[
\pi_t\equiv1\pmod{d_t},\qquad \pi_t\le C_0(4KD)^{L_0}
\]
can be found in time polynomial in the numerical
quantity \(K+D\).
\end{lemma}
\begin{proof}
Let \(H\) be the least power of two at least \(2K\), so
\(2K\le H<4K\). Assign distinct odd residues \(2t-1\pmod H\) to the
indices. Since \(d_t\) is odd, compute
\[
j_t=((2t-2)d_t^{-1})\bmod H\in\{0,\ldots,H-1\},
\qquad c_t=1+d_tj_t.
\]
Then \(1\le c_t<d_tH\), \(c_t\equiv1\pmod{d_t}\), and
\(c_t\equiv2t-1\pmod H\). Thus \(c_t\) is coprime to both \(d_t\)
and \(H\), including when \(d_t\) is composite. By~\eqref{eq:linnik},
the progression \(c_t,c_t+d_tH,c_t+2d_tH,\ldots\) contains a prime
at most \(C_0(d_tH)^{L_0}\). Enumerate it, skipping 1 and testing by
trial division, until the first prime occurs. The bound is polynomial
in \(K+D\), so enumeration, trial division, and the Euclidean inverse
computation all have polynomial running time in that quantity.
The different residues \(2t-1\pmod H\) make the
resulting primes distinct.
\end{proof}

Apply the lemma first to \(s\) inputs, all equal to 37. The number of
inputs is \(s\) and the largest modulus is 37, so its bound gives distinct
variable primes \(q_i\equiv1\pmod{37}\) with
\(q_i\le C_0(148s)^{L_0}\). These primes exceed 37, as required by
Lemmas~\ref{lem:variable} and~\ref{lem:clause}. Since \(s=O(n)\) and
\(C_0,L_0\) are fixed, their numerical values are polynomially bounded
in \(n\).

For the second application, number the \(K=s+m\) tests and set
\(d_t=q_i\) for the test of variable \(\xi_i\), and
\(d_t=q_iq_jq_\ell\) for the test of the clause
\(\xi_i+\xi_j+\xi_\ell=1\). Thus \(\pi_t\equiv1\pmod{d_t}\)
is precisely the field-prime hypothesis of the corresponding test lemma.
Every \(d_t\) is odd and at least three, and
\(D=\max_t d_t\le(\max_i q_i)^3\) is polynomially bounded in \(n\).
The lemma therefore gives distinct field primes
\(\pi_t\equiv1\pmod{d_t}\) with \(\pi_t\le C_0(4KD)^{L_0}\),
again polynomially bounded because \(K=O(n)\).
Both prime searches take polynomial time in \(n\), since their input
counts and largest moduli are polynomially bounded.

\subsection{Construction and correctness of the reduction}
Lemma~\ref{lem:variable} supplies a coefficient congruence for each
variable \(\xi_i\), accepting an exponent \(k\ge0\) precisely when
\(k\bmod q_i\in\{0,1\}\). Under these Boolean restrictions,
Lemma~\ref{lem:clause} supplies a coefficient congruence for each clause,
accepting \(k\) precisely when exactly one of its three residues is one.
Apply the corresponding lemma to each of the \(K\) tests with its
assigned field prime \(\pi_t\). This produces five residues
\((A_t,B_t,a_t,b_t,r_t)\in\F_{\pi_t}^5\) for each test; the numerical
bounds on the primes make all these constructions polynomial in \(n\).
These parameters lie in different finite fields. The following lemma
uses the Chinese remainder theorem to choose one admissible integer
tuple whose coefficient congruence modulo \(M=\prod_t\pi_t\) is
equivalent to simultaneous acceptance by all variable and clause tests.

\begin{lemma}[Combining the local tests]\label{lem:combined-tests}
There exists an admissible tuple \((A,B,a,b,r,M)\) with
\(M=\prod_{t=1}^K\pi_t\), \(A=O(M^2)\), and
\[
(A,B,a,b,r)\equiv(A_t,B_t,a_t,b_t,r_t)\pmod{\pi_t}
\qquad(t=1,\ldots,K).
\]
Let \(U_k,V_k,W_k\) be the coefficient sequence~\eqref{eq:coefficients}
for these parameters \(A,B\). For every \(k\ge0\), the congruence
\(W_k-aU_k-bV_k\equiv r\pmod M\) holds if and only if all local
tests accept \(k\).
\end{lemma}
\begin{proof}
Put
\begin{equation}\label{eq:M}
M=\prod_{t=1}^K\pi_t.
\end{equation}
Since \(K\ge1\), we have \(M\ge2\).
The Chinese remainder theorem gives five representatives
\[
A_0,B,a,b,r\in\{0,\ldots,M-1\}
\]
with the prescribed residues.
Choose
\begin{equation}\label{eq:Alift}
A=A_0+M(M+3).
\end{equation}
This preserves the prescribed residues of \(A\) and gives
\(A\ge(M+1)^2\), because \(M(M+3)-(M+1)^2=M-1\ge1\).
Also \(A<M(M+3)+M=O(M^2)\).
Together with \(0\le B,a,b,r<M\), these bounds show that the tuple
is admissible.

Reduction of the global recurrence~\eqref{eq:recurrence} modulo
\(\pi_t\) gives the corresponding local coefficient sequence. This
follows by induction from the common initial triple \((1,0,0)\) and
the prescribed residues of \(A,B\). By~\eqref{eq:localcongruence}
and the residues of \(a,b,r\), the local test at \(\pi_t\) therefore
accepts \(k\) exactly when
\(W_k-aU_k-bV_k\equiv r\pmod{\pi_t}\).
Since the field primes are distinct, these congruences hold simultaneously
if and only if \(M\) divides \(W_k-aU_k-bV_k-r\).
Thus, for every \(k\ge0\),
\begin{equation}\label{eq:alltests}
W_k-aU_k-bV_k\equiv r\pmod M
\quad\Longleftrightarrow\quad
\text{all local tests accept }k.
\end{equation}
\end{proof}

With these parameters, we output the cubic threshold instance
\((f,\gamma)=(G_*,1)\), where \(G_*\) is the polynomial
in~\eqref{eq:G-star} expanded into monomials. The target instance asks
whether
\[
\min_{(x,y,z)\in\Zp^3}G_*(x,y,z)\le1.
\]

\begin{proof}[Proof of Theorem~\ref{thm:main}]
\emph{Correctness.}
If \(\Phi\) has a satisfying assignment \((\xi_1,\ldots,\xi_s)\), the
Chinese remainder theorem for the distinct variable primes gives
\(k\ge0\) with \(k\equiv\xi_i\pmod{q_i}\) for every \(i\).
All variable tests accept. Every clause test accepts because its three
Boolean residues sum to one. By~\eqref{eq:alltests}, the combined
congruence holds, and Proposition~\ref{prop:equivalence} gives a
nonnegative integer point with \(G_*=1\).

Conversely, if \(G_*\le1\) at a point in \(\Zp^3\), that proposition
supplies an exponent satisfying the combined congruence. By
\eqref{eq:alltests} all variable tests accept, so
\(\xi_i:=k\bmod q_i\in\{0,1\}\). These Boolean residues meet the
hypotheses of every clause test. Each clause therefore has exactly one
true variable. Thus
\begin{equation}\label{eq:reductioniff}
\Phi\text{ is satisfiable}
\quad\Longleftrightarrow\quad
\exists(x,y,z)\in\Zp^3:\ G_*(x,y,z)\le1.
\end{equation}
This argument includes the empty conjunction with its dummy variable.
Proposition~\ref{prop:equivalence} gives positivity on the whole orthant
and exact degree three, as well as a finite, attained minimum.

\emph{Running time and binary size.}
Both prime searches and the constructions of the variable and clause
tests in Lemmas~\ref{lem:variable} and~\ref{lem:clause} take polynomial
time in \(n\), since the number of tests and the numerical values of all
primes are polynomially bounded in \(n\). There are \(K=O(n)\)
field primes, each at most \(n^{O(1)}\), up to an absolute constant, so
\begin{equation}\label{eq:bitM}
\log M=\sum_{t=1}^K\log\pi_t=O(K\log n).
\end{equation}
The standard Chinese remainder algorithm computes \(M\) and the five
representatives in time polynomial in the total binary encoding length
of the moduli and residues. The representatives lie
in \(\{0,\ldots,M-1\}\), and \(A=O(M^2)\), so every parameter has
\(O(\log M)\) bits.
The fixed-degree expression~\eqref{eq:G-star} has at most
\(\binom{6}{3}=20\) monomials. In~\eqref{eq:normpoly}, substituting
\(A=O(M^2)\), \(B,a,b,r=O(M)\) and an affine third coordinate with
coefficients \(O(M)\) bounds every expanded coefficient by \(O(M^5)\).
The expansion therefore has \(O(\log M)\)-bit coefficients and can be
computed in polynomial time.

We have a deterministic polynomial-time many-one reduction, proving the
theorem.
\end{proof}

An exact optimal-value algorithm would distinguish minimum one from
minimum at least two in the constructed instances. Thus computing the
exact optimal value is NP-hard.

\begin{corollary}[NP-completeness of the six-parameter family]\label{thm:NPcomplete}
Let \(\mathcal A_*\) be the language of binary-encoded tuples
\((A,B,a,b,r,M)\) satisfying~\eqref{eq:admissible-star} for which
\(G_*\) takes a value at most one on \(\Zp^3\).
Then \(\mathcal A_*\) is NP-complete. Every yes-instance has an exponent
certificate \(0\le k<M^3\), verifiable in time polynomial in the tuple's
binary encoding length.
\end{corollary}
\begin{proof}
NP-hardness follows from the tuple construction and~\eqref{eq:reductioniff}.
By Proposition~\ref{prop:equivalence}, an admissible tuple belongs to
\(\mathcal A_*\) exactly when some exponent \(k\ge0\) satisfies
\(W_k-aU_k-bV_k\equiv r\pmod M\). The periodicity argument in the proof
of Proposition~\ref{prop:equivalence} gives a period \(h\le M^3\) for the
coefficient sequence modulo \(M\). Reducing an
accepted exponent modulo \(h\) preserves the congruence and gives
\(0\le k<h\le M^3\), using at most \(3\lceil\log_2 M\rceil\) bits.

Since admissibility is part of the definition of \(\mathcal A_*\), a
verifier first rejects tuples that fail~\eqref{eq:admissible-star} and
certificates outside \(0\le k<M^3\). It then computes
\(C^k(1,0,0)^{\mathsf T}\pmod M\) by repeated squaring, with \(C\)
as in~\eqref{eq:matrix}. By~\eqref{eq:recurrence}, this gives
\((U_k,V_k,W_k)^{\mathsf T}\) modulo \(M\), so the verifier can test the
congruence. Proposition~\ref{prop:equivalence} ensures that every accepted
certificate yields a nonnegative integer point with \(G_*=1\).
The exponentiation uses \(O(\log M)\) multiplications of \(3\times3\)
matrices on \(O(\log M)\)-bit residues. Checking admissibility and
initially reducing \(A\) modulo \(M\) take time polynomial in the tuple's
binary encoding length. Thus verification is polynomial-time, proving
membership in NP.
\end{proof}

\section{Conclusion}\label{sec:conclusion}
Our construction proves NP-hardness of cubic integer optimization over
the fixed orthant \(\Zp^3\), even at threshold one for objectives with
positive integer values and a finite, attained minimum. The proof first
constructs a family of cubic objectives from a norm form and shows that
each objective has minimum one exactly when an associated coefficient
congruence has a solution. Finite-field tests then encode Boolean
variables and exactly-one clauses as congruences of this form. The Chinese
remainder theorem combines these tests and determines the parameters of
the objective used in the reduction.

For the six-parameter family constructed here, every yes-instance admits
a short certificate. A remaining challenge is to obtain polynomial-size,
polynomial-time verifiable certificates for arbitrary cubic objectives
over \(\Zp^3\). We conjecture that the language \(\mathcal L_3\)
in~\eqref{eq:threshold} belongs to NP. An even more challenging open
problem is whether the corresponding threshold problem over \(\Zp^n\)
belongs to NP when \(n\) is part of the input.

\section*{Statements and Declarations}
\paragraph{Funding.}
Linchuan Wei was supported by the University Development Fund of The Chinese
University of Hong Kong, Shenzhen, under Grant No. UDF01003946.

\paragraph{Competing interests.}
The author has no competing interests to declare.

\paragraph{Data availability.}
No datasets were generated or analysed during this study.

\paragraph{Use of AI.}
The authors used GPT to assist with reviewing the manuscript, refining its
exposition, and preparing the bibliography. They reviewed all AI-assisted
material and take full responsibility for the manuscript.

\bibliographystyle{spmpsci}
\bibliography{references}
\end{document}